\documentclass[cmbright,doublespace]{mmaauth}
\usepackage{moreverb}
\usepackage[dvips,colorlinks,bookmarksopen,bookmarksnumbered,citecolor=red,urlcolor=red]{hyperref}

\newcommand\BibTeX{{\rmfamily B\kern-.05em \textsc{i\kern-.025em b}\kern-.08em
T\kern-.1667em\lower.7ex\hbox{E}\kern-.125emX}}

\def\volumeyear{2016}

\DeclareMathOperator{\real}{Re}
\newcommand{\re}{\mathrm{e}}
\newcommand{\ri}{\mathrm{i}}
\newcommand{\rd}{\mathrm{d}}
\newcommand{\DD}[2]{\frac{\partial #1}{\partial #2}}

\newcommand{\DDM}[3]{\frac{\partial^2 #1}{{\partial #2}{\partial #3}}}
\newcommand*{\doi}[1]{\href{http://dx.doi.org/\detokenize{#1}}{doi:\detokenize{#1}}}
\allowdisplaybreaks

\begin{document}

\runninghead{I.~C.~Christov and C.~I.~Christov}

\title{On mechanical waves and Doppler shifts from moving boundaries}

\author{Ivan~C.~Christov\affil{a}\footnotemark[2]\corrauth\ and C.~I.~Christov\affil{b}\footnotemark[3]} 
\address{\affilnum{a}Theoretical Division and Center for Nonlinear Studies,
Los Alamos National Laboratory, Los Alamos, New Mexico 87545, USA\\
\affilnum{b}Department of Mathematics, University of Louisiana at Lafayette, Lafayette, Louisiana 70504, USA}
\corraddr{School of Mechanical Engineering, Purdue University, West Lafayette, Indiana 47907, USA. E-mail: \href{mailto:christov@purdue.edu}{christov@purdue.edu}; URL: \url{http://christov.tmnt-lab.org}.}
\cgs{Los Alamos National Laboratory (LANL) is operated by Los Alamos National Security, L.L.C.\ for the National Nuclear Security Administration of the U.S.\ Department of Energy under Contract No.\ DE-AC52-06NA25396.}

\begin{abstract}
We investigate the propagation of infinitesimal harmonic mechanical waves emitted from a boundary with variable velocity and arriving at a stationary observer. In the classical Doppler effect,  $X_\mathrm{s}(t) = vt$ is the location of the source with constant velocity $v$. In the present work, however, we consider a source co-located with a moving boundary $x=X_\mathrm{s}(t)$, where $X_\mathrm{s}(t)$ can have an arbitrary functional form. For ``slowly moving'' boundaries (\textit{i.e.}, ones for which the timescale set by the mechanical motion is large in comparison to the inverse of the frequency of the emitted wave), we present a multiple-scale asymptotic analysis of the moving-boundary problem for the linear wave equation. We obtain a closed-form leading-order (with respect to the latter small parameter) solution and show that the variable velocity of the boundary results not only in frequency modulation but also in amplitude modulation of the received signal. Consequently, our results extending the applicability of two basic tenets of the theory of a moving source on a stationary domain, specifically that (a) $\dot X_\mathrm{s}$ for non-uniform boundary motion can be inserted in place of the constant velocity $v$ in the classical Doppler formula and (b) that the non-uniform boundary motion introduces variability in the amplitude of the wave. The specific examples of decelerating and oscillatory boundary motion are worked out and illustrated.
\end{abstract}

\MOS{35L05, 35C20, 78M35}

\keywords{Doppler effect; accelerating source; multiple-scales expansion; wave equation; moving boundary}

\maketitle

\footnotetext[2]{Present address: School of Mechanical Engineering, Purdue University, West Lafayette, Indiana 47907, USA.}
\footnotetext[3]{Prof.~C.~I.~Christov passed away prior to submission of the manuscript. This paper is dedicated to his memory.}

\vspace{-6pt}

\section{Introduction}
\label{sec:intro}

The classical (or non-relativistic) \emph{Doppler effect} \cite{doppler} is concerned with the change in observed frequency of a mechanical wave when its emitter is in relative motion with respect to the observer \cite{g65,eh69}. The shifted frequency $\omega_\mathrm{D}$, measured by a stationary observer downstream (along the direction of propagation of the wave) from the emitter, goes (see, \textit{e.g.}, \cite{g65}) as
\begin{equation}
\omega_\mathrm{D} = \frac{\omega}{1-v/c},
\label{eq:doppler_shift}
\end{equation}
where $\omega$ is the frequency of the emitted waves, $v$ is the velocity of the source, and $c$ is the phase speed of infinitesimal waves in the particular medium under consideration (\textit{e.g.}, the speed of sound for an acoustic wave). The quantity $\Delta\omega := \omega_\mathrm{D} - \omega$ is termed the \emph{Doppler shift}. For a source moving towards the observer ($v>0$), the measured frequency $\omega_\mathrm{D}$ is larger than the emitted frequency $\omega$, whereas for a source moving away from the observer ($v<0$), $\omega_\mathrm{D} < \omega$. The Doppler effect is a staple of classical wave physics, and the applications of \eqref{eq:doppler_shift} in telecommunications, meteorology, medicine, etc.\ are so numerous that we do not attempt to list them here (see, \textit{e.g.}, \cite{g65,eh69,mi68,w74}). 
However, more than 150 years after Christian Doppler's proposal, novel aspects of the Doppler effect continue to be uncovered \cite{p14}. 

The case of translation of the source along the line $x=X_\mathrm{s}(t)\equiv vt$ at constant velocity $v$ has been exhaustively treated and is well understood. For non-uniform source velocities, \textit{i.e.}, $\dot X_\mathrm{s}(t) \ne const.$, the problem has been analyzed more recently \cite{o97,ow16} due to its relevance for acoustics in moving and inhomogeneous media \cite{o85,o89,o94,oetal04}. Simple acoustic laboratory experiments \cite{bf01,a07} have been performed showing the effects of acceleration of the emitter. Specifically, higher harmonics appear in the spectrum of the received signal, and the spectrum itself becomes markedly asymmetric. In the analysis of the data, however, it is common to formally replace $v$ with $\dot X_\mathrm{s}(t)$ in \eqref{eq:doppler_shift} \cite{bf01,a07}. This manipulation is justified by the solutions, given in \cite[Chap.~5]{ow16}, of the three-dimensional (3D0 initial-boundary-value problem (IBVP) for the wave equation with the acoustic source modeled as a singular term on the right hand side of the linear wave equation \cite{o97,ow16,gvn09}. Here, we would like to pose a different variant of this problem: if the source were attached to a \emph{moving} boundary of a one-dimensional (1D) domain, then can we still replace $v$ with $\dot X_\mathrm{s}(t)$ in \eqref{eq:doppler_shift}, or would corrections arise from a formal mathematical analysis? Posing the problem in this manner also provides a natural generalization of some fundamental IBVPs, which we review in context below, studied in the mathematics literature.

From the mathematical point of view, some basic sketches of the theory of such mechanical wave motions have been presented in the context of asymptotic and perturbation methods \cite{n73,m81,h95,kc96}. However, to the best of our knowledge, an analysis of wave propagation from an emitter co-located with a boundary of non-uniform velocity cannot be found in the literature. Thus, in this paper, we provide a formal perturbative solution based on the method of \emph{multiple scales} for the case in which the timescale set by the mechanical motion of the emitter is large compared to the inverse frequency of the emitted wave. 
We show that for general boundary velocity given by $\dot{X}_\mathrm{s}(t)$, the expression for the shifted frequency \eqref{eq:doppler_shift} can be immediately modified as
\begin{equation}
\omega_\mathrm{D}(t) = \frac{\omega}{1-\dot{X}_\mathrm{s}(t)/c},
\label{eq:doppler_shift_acc_heur}
\end{equation}
within the assumed order of approximation. This result is, of course, exactly in agreement with the corresponding Doppler shift found from the general 3D solution for a moving source in a homogeneous stationary medium \cite[eq.~(5.23)]{ow16}.
An additional physical effect obtained by our analysis is that the amplitude of the wave is also affected by the non-uniform boundary motion, which can also be inferred from the the general 3D solution for a moving source in a homogeneous stationary medium \cite[\S5.1]{ow16}.

To this end, in \S\ref{sec:position}, we reformulate the moving-boundary problem for the wave equation into an equivalent problem on a fixed domain for a dispersive, variable-coefficient wave equation. In \S\ref{sec:ms} we give the leading-order solution by the method of multiple scales. Then, in \S\ref{sec:discuss}, the solution is illustrated for a decelerating boundary and for an oscillatory boundary motion. Finally, in \S\ref{sec:conclusion}, conclusions are stated and a broader context for the present results and their applicability is proposed.

\section{Position of the problem}
\label{sec:position}

A plethora of mechanical wave phenomena are governed by the $(1+1)$ dimensional linear wave partial differential equation (PDE)
\begin{equation}
\frac{\partial^2 U}{\partial t^2} - c^2 \frac{\partial^2 U}{\partial x^2} = 0,\qquad X_\mathrm{s}(t) < x < \infty,\quad 0 < t < \infty,
\label{eq:wave_eq}
\end{equation}
where $U=U(x,t)$ can be, \textit{e.g.}, 
the acoustic potential \cite{mi68}, the elastic displacement \cite{a73} or even the temperature field (under certain nonclassical theories of thermoelasticity) \cite[\S2.3]{s11}, to name a few, and $c$ is the phase speed of infinitesimal waves in the material medium. In the present work, \eqref{eq:wave_eq} is subject to the boundary condition
\begin{equation}
U\big(X_\mathrm{s}(t),t\big) = \re^{\ri\omega t},\qquad X_\mathrm{s}(t) := vt + a\alpha(\Omega t),\qquad t>0,
\label{eq:wave_eq_bc}
\end{equation}
which represents an \emph{accelerating} moving source emitting monochromatic harmonic waves with frequency $\omega$, where $v$, $a$ and $\Omega$ are some positive constants and $\alpha$ is a dimensionless function (nonlinear in its argument). Here, $a$ has the dimension of length and $\Omega$ has the dimension of inverse time. We have chosen this particular functional form for $X_\mathrm{s}(t)$ so that the classical Doppler effect is easily recovered in the final results by setting $\alpha = 0$ and/or taking the limit $a\to0^+$. We take $\alpha$ such that $\alpha(0)=0$, without loss of generality, so that $X_\mathrm{s}(0)=0$. By convention, we work with complex exponentials since \eqref{eq:wave_eq} is a linear equation, and the real part of $U$ is taken at the end of the calculation.

In addition, we must supplement eqs.~\eqref{eq:wave_eq} and \eqref{eq:wave_eq_bc} with the radiation condition
\begin{equation}
\lim_{x\to\infty} \left(\frac{\partial U}{\partial x} + \ri \kappa U \right) = 0,
\label{eq:wave_eq_rad_bc}
\end{equation}
where $\kappa$ is the spatial wave number, and the ``$+$'' is chosen so that only waves that are outgoing at $x=\infty$ are allowed \cite[\S28]{s49}. Without loss of generality, homogeneous initial conditions, $U(x,0)=\frac{\partial U}{\partial t}(x,0)=0$, can be imposed because, for the present purposes, we are only interested in the influence of the boundary condition, meaning that $\forall x<\infty$ $\exists t = x/c$ such that any non-zero initial condition has propagated past this location, and only the effects due to the boundary condition are ``felt'' there.\footnote{This type of argument could be generalized to ``arbitrary'' initial conditions as long as they are constrained to produce waves satisfying the radiation condition \eqref{eq:wave_eq_rad_bc}.}

George Carrier's ``spaghetti problem'' \cite{c49} regarding the normal modes of a string being shortened due to its accelerated withdrawal into an orifice motivated some early analytical work by Balazs \cite{b61} and Greenspan \cite{g63} on the Dirichlet IBVP for \eqref{eq:wave_eq} on a finite domain with moving and/or accelerating boundaries. A pernicious  feature of these problems is reflections from the boundaries, leading to analytical solutions in the form of trigonometric series \cite{b61,g63,c93}. The Dirichlet problem can also be solved for general boundary motions using nonlinear transformations of the independent variables \cite{g98,wh89} and multiple-scale asymptotics \cite{gm81}. Integral representations for the solution to the half-space problem with a moving boundary have also recently been proposed \cite{pp10} on the basis of advanced transform techniques \cite{fok97,fp00,fok03,pel06}. In contrast, we study the physical (rather than abstract) half-space problem for harmonic mechanical waves in order to discern any frequency and/or amplitude shifts due to the non-uniform motion of the boundary.

We choose to convert eqs.~\eqref{eq:wave_eq}--\eqref{eq:wave_eq_rad_bc} into a boundary-value problem on $[0,\infty)$ by introducing the moving frame coordinate 
\begin{equation}
\xi = x - X_\mathrm{s}(t) \equiv x - vt - a\alpha(\Omega t),
\end{equation}
while keeping the time coordinate the same. Then, the temporal and spatial partial derivatives transform as
\begin{equation}
\frac{\partial^2}{\partial t^2} = \frac{\partial^2}{\partial t^2} - 2(v+a\Omega \alpha')\frac{\partial^2}{\partial\xi\partial t} - a\Omega^2\alpha''\frac{\partial}{\partial \xi} + (v+a\Omega \alpha')^2\frac{\partial^2}{\partial \xi^2},\qquad \frac{\partial^2 }{\partial x^2} = \frac{\partial^2 }{\partial \xi^2},
\label{eq:pdiff_moving}
\end{equation}
where a prime indicates differentiation with respect to \emph{the argument} of $\alpha$.
Letting $U(x,t) = \tilde{U}(\xi,t)$ and introducing \eqref{eq:pdiff_moving} into \eqref{eq:wave_eq}, we obtain
\begin{equation}
\frac{\partial^2 \tilde{U}}{\partial t^2} - 2(v+a\Omega\alpha')\frac{\partial^2 \tilde{U}}{\partial\xi\partial t} - a\Omega^2\alpha'' \frac{\partial \tilde{U}}{\partial \xi} - \big[c^2 - (v+a\Omega\alpha')^2 \big]\frac{\partial^2 \tilde{U}}{\partial \xi^2} = 0,\qquad
 0<\xi<\infty,\quad 0<t<\infty.
\label{eq:wave_eq_mf}
\end{equation}
Equation~\eqref{eq:wave_eq_mf} is now a \emph{dispersive} wave equation with \emph{variable} coefficients. The general theory of dispersive waves under such equations is described by Whitham \cite[Chap.~11]{w74}. Some  remarks on the theory of such PDEs, including analysis of the Lie symmetries, were given by Bluman \cite{b83}. Analytical solutions for special choices of the coefficient have been provided, \textit{e.g.}, for $v=a=0$ and $c = c(x)$ \cite{bk87,bk88,t90} or $v=a=0$ and $c = c(x,t)$ \cite{SV87}, with further generalization given in \cite{SV88}. Others have considered the case of discontinuous $c(x)$ \cite{d72}. More recently, the asymptotic properties of \emph{localized} solutions for $v=a=0$ but $c = c(x)$ have been examined in detail \cite{dst07,dst08}. Some special constant coefficient cases of \eqref{eq:wave_eq_mf} arise in the study of low-frequency modulation of acoustic radiation forces \cite{dg11}. Unfortunately, all these results are too specialized to be of immediate use in our analysis of \eqref{eq:wave_eq_mf}.

Before proceeding further, we must enforce some limitations on $v$, $a$ and $\Omega$. Obviously, at any time $t$, we must have 
\begin{equation}
|\dot X_\mathrm{s}(t)| \equiv |v + a\Omega\alpha'(\Omega t)| < c, 
\label{eq:sound_speed_ineq}
\end{equation}
\textit{i.e.}, the instantaneous velocity of the boundary must be less that the phase speed of waves, otherwise \eqref{eq:wave_eq_mf} changes type from hyperbolic to elliptic, and the problem becomes ill-posed (and unphysical).  In other words, the boundary motion is \emph{subsonic}.\footnote{Note that, in the case of the IBVP for wave propagation in a homogeneous medium at rest, in which the acoustic source modeled as a singular term on the right hand side of \eqref{eq:wave_eq} (see, \textit{e.g.}, \cite[Chap.~5]{ow16}), it possible to also consider \emph{supersonic} sources.} We consider the case when $a\Omega/c$ (a type of acceleration-based Strouhal number) is $\mathcal{O}(1)$, \textit{i.e.}, the time scale set by the acceleration of the source and the time scale on which its acceleration varies are comparable. Since we introduced the parameter $a$, we are free to normalize $\alpha$ so that $\max_{\mathfrak{t}\ge0}|\alpha'(\mathfrak{t})| = 1$. Consequently, a necessary condition for the  inequality \eqref{eq:sound_speed_ineq} to hold is $a\Omega/c < 1 - v/c$.

The most important assumption we make, however, is that $\omega \gg \Omega$, \textit{i.e.}, the frequency of the emitted wave $\omega$ is much larger than the frequency of the mechanical motion $\Omega$ associated with the acceleration of the source.\footnote{Note that, in the case of the IBVP for wave propagation in a homogeneous medium at rest, in which the acoustic source modeled as a singular term on the right hand side of \eqref{eq:wave_eq} (see, \textit{e.g.}, \cite[Chap.~5]{ow16}), it possible to find an analytical solution without requiring that $\omega \gg \Omega$.} (Equivalently, the timescale $\Omega^{-1}$ set by the mechanical motion is large in comparison to the time scale $\omega^{-1}$ of wave propagation; \textit{i.e.}, a ``slowly moving'' source.) This assumption defines the small parameter for the upcoming asymptotic expansion. Therefore, we introduce the following dimensionless independent variables and dimensionless parameters: 
\begin{equation}
\tau = \omega t,\qquad \eta = \xi \omega/c,\qquad \beta := v/c,\qquad \delta := a\Omega/c,\qquad \epsilon := {\Omega}/{\omega} \ll 1.
\label{eq:var_param}
\end{equation}
Note that the non-dimensionalization of the dependent variable is arbitrary since \eqref{eq:wave_eq_mf} is a homogeneous linear equation, hence it is invariant under re-scaling of the dependent variable. 

Now, letting $\tilde{U}(\xi,t) = \hat{U}(\eta,\tau)$ and making use of these dimensionless variables from \eqref{eq:var_param}, \eqref{eq:wave_eq_mf} becomes
\begin{equation}
\frac{\partial^2 \hat{U}}{\partial \tau^2} - 2 \big[\beta +  \delta\alpha'(\epsilon\tau) \big] \frac{\partial^2 \hat{U}}{\partial\eta\partial\tau} - \epsilon\delta\alpha''(\epsilon\tau)\frac{\partial\hat{U}}{\partial\eta} - \left\{ 1 - \big[\beta + \delta\alpha'(\epsilon\tau) \big]^2 \right\}\frac{\partial^2 \hat{U}}{\partial \eta^2} = 0,\qquad
 0<\eta<\infty,\quad 0<\tau<\infty.
\label{eq:wave_eq_nd}
\end{equation}
Recall that, here, primes stand for differentiation of a function with respect to \emph{its argument} (in this case, $\epsilon\tau$).

Finally, the boundary and radiation condition from eqs.~\eqref{eq:wave_eq_bc} and \eqref{eq:wave_eq_rad_bc} become
\begin{equation}
\hat{U}(0,\tau) = \re^{\ri \tau},\qquad  \lim_{\eta\to\infty} \left(\frac{\partial \hat{U}}{\partial \eta} + \ri k \hat{U} \right) = 0,
\label{eq:wave_eq_bc_nd}
\end{equation}
where $k = \kappa c/\omega$ is the dimensionless wave number.

Hence, we have transformed our original movin- boundary problem into a variable-coefficient problem on the half-line. Unfortunately, our problem does not appear to yield itself to a closed-form solution. However, the variable coefficients in \eqref{eq:wave_eq_nd} are \emph{slowly} varying, \textit{i.e.}, they depend only on $\epsilon \tau$. Thus, we proceed by perturbation methods as in \cite{gm81,dg11,kk59}.

\section{Solution by a multiple-scales expansion}
\label{sec:ms}

Equation~\eqref{eq:wave_eq_nd} is a linear wave equation with slowly varying coefficients, which makes it an ideal candidate for a multiple-scales asymptotic expansion \cite{n73,m81,h95,kc96,e00,gk09}, the generalization of Cole's two-variable expansion procedure \cite[Chap.~3]{Cole}. To this end, we introduce the ``fast'' time $t_0 = \tau$, the ``slow'' time $t_1 = \epsilon \tau$, the ``short'' spatial coordinate $y_0 = \eta$ and the ``long'' spatial coordinate $y_1 = \epsilon \eta$. For convenience, we first rewrite \eqref{eq:wave_eq_nd} as
\begin{equation}
\frac{\partial^2 \hat{U}}{\partial \tau^2} - 2\tilde\beta(t_1)\frac{\partial^2 \hat{U}}{\partial\eta\partial\tau} 
 - \big[1 - \tilde\beta^2(t_1) \big] \frac{\partial^2 \hat{U}}{\partial \eta^2} 
 - \epsilon\delta\alpha''(t_1)\frac{\partial\hat{U}}{\partial\eta} = 0,
\label{eq:wave_eq_nd_1}
\end{equation}
where $\tilde\beta(t_1) := \beta + \delta\alpha'(t_1)$ is a function of the slow time alone. Then, we let $\hat{U}(\eta,\tau) = \mathcal{U}(y_0,y_1,t_0,t_1)$ with its partial derivatives transforming as
\begin{equation}
\begin{aligned}
\frac{\partial \hat{U}}{\partial \tau} &= \frac{\partial \mathcal{U}}{\partial t_0} + \epsilon \frac{\partial \mathcal{U}}{\partial t_1},\qquad
\frac{\partial^2 \hat{U}}{\partial \tau^2} = \frac{\partial^2 \mathcal{U}}{\partial t_0^2} + 2\epsilon \frac{\partial^2 \mathcal{U}}{\partial t_0\partial t_1} + \epsilon^2 \frac{\partial^2 \mathcal{U}}{\partial t_1^2},\\ 
\frac{\partial \hat{U}}{\partial \eta} &= \frac{\partial \mathcal{U}}{\partial y_0} + \epsilon \frac{\partial \mathcal{U}}{\partial y_1},\qquad
\frac{\partial^2 \hat{U}}{\partial \eta^2} = \frac{\partial^2 \mathcal{U}}{\partial y_0^2} + 2\epsilon \frac{\partial^2 \mathcal{U}}{\partial y_0\partial y_1} + \epsilon^2 \frac{\partial^2 \mathcal{U}}{\partial y_1^2},\\
\frac{\partial^2 \hat{U}}{\partial \eta \partial \tau} &= \frac{\partial^2 \mathcal{U}}{\partial y_0\partial t_0} + \epsilon \frac{\partial^2 \mathcal{U}}{\partial y_0 \partial t_1} + \epsilon \frac{\partial^2 \mathcal{U}}{\partial y_1\partial t_0} + \epsilon^2 \frac{\partial^2 \mathcal{U}}{\partial y_1 \partial t_1}.
\end{aligned}
\label{eq:twotime_twospace_deriv}
\end{equation}
Upon introducing \eqref{eq:twotime_twospace_deriv} into \eqref{eq:wave_eq_nd_1} and keeping only leading-order terms and terms proportional to $\epsilon$, we obtain
\begin{multline}
\frac{\partial^2 \mathcal{U}}{\partial t_0^2} 
- 2\tilde\beta(t_1)\frac{\partial^2 \mathcal{U}}{\partial y_0 \partial t_0} 
- \big[1-\tilde\beta^2(t_1) \big]\frac{\partial^2 \mathcal{U}}{\partial y_0^2}
+ 2\epsilon \frac{\partial^2 \mathcal{U}}{\partial t_0\partial t_1} 
- 2\epsilon\tilde\beta(t_1)\frac{\partial^2 \mathcal{U}}{\partial y_0 \partial t_1}\\
- 2\epsilon\tilde\beta(t_1)\frac{\partial^2 \mathcal{U}}{\partial y_1 \partial t_0}
- 2\epsilon \big[1-\tilde\beta^2(t_1) \big]\frac{\partial^2 \mathcal{U}}{\partial y_0 \partial y_1} 
- \epsilon\delta\alpha''(t_1)\frac{\partial\mathcal{U}}{\partial y_0} + \mathcal{O}(\epsilon^2) = 0.
\label{eq:wave_eq_ms_vars1}
\end{multline}

We proceed in the usual manner by a regular expansion of the dependent variable:
\begin{equation}
\mathcal{U}(y_0,y_1,t_0,t_1) = \mathcal{U}_0(y_0,y_1,t_0,t_1) + \epsilon \mathcal{U}_1(y_0,y_1,t_0,t_1) + \cdots.
\end{equation}
In turn, the (first) boundary condition from \eqref{eq:wave_eq_bc_nd} becomes
\begin{equation}
\mathcal{U}_0(0,0,t_0,t_1) = \re^{\ri t_0},\qquad \mathcal{U}_j(0,0,t_0,t_1) = 0 \quad (j>0).
\label{eq:ms_bc}
\end{equation}
Then, at the leading order, \eqref{eq:wave_eq_ms_vars1} becomes
\begin{equation}
\mathcal{L}_0[\mathcal{U}_0] = 0,\qquad \mathcal{L}_0 := \frac{\partial^2}{\partial t_0^2} - 2\tilde\beta\frac{\partial^2}{\partial y_0 \partial t_0} - (1-\tilde\beta^2)\frac{\partial^2}{\partial y_0^2}.
\label{eq:ms_ord1_pb}
\end{equation}
Here, it is important to recall that $\tilde\beta$ does \emph{not} depend on the fast time $t_0$. Therefore, according to the multiple-scales expansion procedure, it is considered a constant at this order. This is equivalent to the assumption that the wavenumber $k$ depends on the slow time $t_1$, which is sometimes referred to as the ``generalized'' method of multiple scales \cite[\S6.4]{n73}.

Clearly, a solution of eqs.~\eqref{eq:ms_ord1_pb} and \eqref{eq:ms_bc} of the form\footnote{Here, we use the structure of the boundary condition to infer the form of the solution. More generally, for arbitrary  excitations, one can introduce characteristic coordinates and proceed along the lines of the Appendix.}
\begin{equation}
\mathcal{U}_0(y_0,y_1,t_0,t_1) = \mathcal{A}_0(y_1,t_1)\re^{\ri t_0 - \ri k y_0 + \ri\psi_0(y_1,t_1)}
\label{eq:U0_solution}
\end{equation}
exists provided that 
\begin{equation}
-1 -2\tilde\beta k + \big(1-\tilde\beta^2 \big)k^2=0 \quad \Longrightarrow \quad k = \frac{\tilde\beta \pm 1}  {1 - \tilde\beta^2} = \frac{\pm 1}{1 \mp \tilde\beta}. 
\label{eq:k_beta}
\end{equation}
Here, we must pick the upper sign in the expression for $k$ ($\Rightarrow k = 1/(1 -\tilde\beta)>0$) to satisfy the radiation condition \eqref{eq:wave_eq_bc_nd} (\textit{i.e.}, to have waves that are propagating away from the source and outgoing at $\eta = \infty$) and set $\mathcal{A}_0(0,t_1)=1$, $\psi_0(0,t_1) = 0$ to satisfy the boundary condition \eqref{eq:ms_bc}. Notice that the earlier assumptions that $a\Omega/c < 1 - v/c\Leftrightarrow \delta < 1-\beta$ and that $\max_{\mathfrak{t}\ge0}|\alpha'(\mathfrak{t})|=1$ guarantees that $(1-\tilde\beta) >0 \Rightarrow 1/(1-\tilde\beta) >0$ in the accelerating case (\textit{i.e.}, the case of $\delta\alpha' > 0$). In general, the sign of $\beta$ depends on whether the boundary is moving towards or away from the observer, the latter being situated somewhere on the positive abscissa. For $\beta >0$, we have $k>1$ (the wave is ``shorter''), \textit{i.e.}, the observed pitch is higher for a source moving towards the observer. 
Conversely, for $\beta <0$, $k<1$ (the wave is ``longer''), which means the observed pitch is lower for a receding source.

Continuing to $\mathcal{O}(\epsilon)$, we must now solve
\begin{equation}
\mathcal{L}_0[\mathcal{U}_1] = 
- 2\frac{\partial^2 \mathcal{U}_0}{\partial t_0\partial t_1} 
+ 2\tilde\beta\frac{\partial^2 \mathcal{U}_0}{\partial y_0\partial t_1} 
+ 2\tilde\beta\frac{\partial^2 \mathcal{U}_0}{\partial y_1\partial t_0} 
+ 2\big(1-\tilde\beta^2 \big) \frac{\partial^2 \mathcal{U}_0}{\partial y_0 \partial y_1}
+ \delta\alpha''\frac{\partial \mathcal{U}_0}{\partial y_0}.
\label{eq:ord_eps_eq}
\end{equation}
Denoting the right-hand side above as $\mathfrak{F}$, it can be evaluated based on the solution for $\mathcal{U}_0$ from \eqref{eq:U0_solution}:
\begin{multline}
\mathfrak{F}(y_0,y_1,t_0,t_1) = \Bigg\{-2\left( \ri \frac{\partial \mathcal{A}_0}{\partial t_1} - \mathcal{A}_0 \frac{\partial \psi_0}{\partial t_1} \right) 
+ 2\tilde\beta\left(-\ri k\frac{\partial \mathcal{A}_0}{\partial t_1} + k \mathcal{A}_0\frac{\partial \psi_0}{\partial t_1} \right)\\
+ 2\tilde\beta\left(\ri \frac{\partial \mathcal{A}_0}{\partial y_1} - \mathcal{A}_0\frac{\partial \psi_0}{\partial y_1} \right)
+ 2\big(1-\tilde\beta^2 \big) \left(-\ri k\frac{\partial \mathcal{A}_0}{\partial y_1} + k \mathcal{A}_0\frac{\partial \psi_0}{\partial y_1} \right)
- \ri \delta\alpha'' k \mathcal{A}_0
\Bigg\} \re^{\ri t_0 - \ri k y_0 + \ri\psi_0}.
\label{eq:ms_ord1_rhs}
\end{multline}
This right-hand side of \eqref{eq:ord_eps_eq} will produce secular terms because $\re^{\ri t_0 - \ri k y_0 + \ri\psi_0}$ is in the nullspace of $\mathcal{L}_0$. Therefore, we must choose $\mathcal{A}_0$ and $\psi_0$ so that $\mathfrak{F}\equiv0$. Separating the real and imaginary parts of \eqref{eq:ms_ord1_rhs} and assuming a nontrivial solution $\mathcal{A}_0\ne 0$, we obtain
\begin{subequations}\begin{align}
(1+\tilde\beta k)\frac{\partial \psi_0}{\partial t_1} - \big[\tilde\beta - \big(1-\tilde\beta^2 \big)k \big] \frac{\partial \psi_0}{\partial y_1}  &= 0,\\
(1+\tilde\beta k)\frac{\partial \mathcal{A}_0}{\partial t_1} - \big[\tilde\beta - \big(1-\tilde\beta^2 \big)k \big] \frac{\partial \mathcal{A}_0}{\partial y_1} &= -\tfrac{1}{2}\delta\alpha'' k \mathcal{A}_0.
\end{align}\label{eq:amp_phase_transport}\end{subequations}

We proceed by the method of characteristics \cite[Chap.~II]{ch62}, which transforms \eqref{eq:amp_phase_transport} into a set of ordinary differential equations (ODEs):
\begin{subequations}\begin{align}
\frac{\rd t_1}{\rd s} &= 1+\tilde\beta k,\\ 
\frac{\rd y_1}{\rd s} &= -\big[\tilde\beta - \big(1-\tilde\beta^2 \big)k \big],\\ 
\frac{\rd \overline{\psi}_0}{\rd s} &= 0,\label{eq:char_odes_psi0}\\
\frac{\rd \mathcal{\overline A}_0}{\rd s} &= - \tfrac{1}{2}\delta \alpha'' k\mathcal{\overline A}_0,\label{eq:char_odes_A0}
\end{align}\label{eq:char_odes}\end{subequations}
where we have introduced the notation $\psi_0(y_1,t_1) = \overline{\psi}_0(s)$ and $\mathcal{A}_0(y_1,t_1) = \mathcal{\overline A}_0(s)$. To find the characteristics, we need to integrate the first two ODEs in \eqref{eq:char_odes}, keeping in mind that $\tilde\beta=\tilde\beta(t_1)$, subject to the ``initial'' condition that $s=0$ when $t_1=0$. Therefore, upon using the expression for $k$ from \eqref{eq:k_beta}, we have
\begin{subequations}\begin{alignat}{3}
\frac{\rd s}{\rd t_1} &= 1 - \tilde\beta{(t_1)} \qquad &&\Longrightarrow \qquad  s(y_1,t_1) &&= \int_0^{t_1} {1- \tilde\beta{(\mathfrak{t})}} \,\rd \mathfrak{t} + \mathfrak{C}_1(y_1)
\label{eq:ds_dt1},\\
\frac{\rd s}{\rd y_1} &= 1 \qquad &&\Longrightarrow \qquad  s(y_1,t_1) &&= y_1 + \mathfrak{C}_2(t_1).
\end{alignat}\end{subequations}
Solving for the arbitrary functions $\mathfrak{C}_{1,2}$ between the two equations and recalling that $\tilde\beta(t_1) \equiv \beta + \delta \alpha'(t_1)$, we obtain
\begin{equation}
s(y_1,t_1) = \int_0^{t_1} {1- \tilde\beta{(\mathfrak{t})}} \,\rd \mathfrak{t} + y_1 = 
(1-\beta)t_1 - \delta\alpha(t_1) + y_1.
\label{eq:s_t1_y1}
\end{equation}
From \eqref{eq:char_odes_A0}, we obtain
\begin{equation}
\ln\big|\overline{\mathcal{A}}_0\big(s(y_1,t_1)\big) \big| - \ln \big |\overline{\mathcal{A}}_0\big(s(0,t_1)\big) \big| = - \frac{1}{2}\int_{s(0,t_1)}^{s(y_1,t_1)}\frac{ \delta\alpha'' (t_1(\varsigma))}{1 - \tilde\beta(t_1(\varsigma))} \,\rd \varsigma,
\end{equation}
where the limits of integration were determined by the fact that our boundary condition is given at $y_0=y_1=0$.
Making a change of variables in the integral by using \eqref{eq:ds_dt1} and recalling that $\overline{\mathcal{A}}_0\big(s(0,t_1)\big) = 1$ from the boundary condition \eqref{eq:ms_bc}, we have
\begin{equation}
\ln|\mathcal{A}_0(y_1,t_1)| = - \frac{1}{2}\int_{s(0,t_1)}^{s(y_1,t_1)} \delta\alpha'' (\mathfrak{t}) \,\rd \mathfrak{t}\quad
\Longrightarrow \quad 
\mathcal{A}_0(y_1,t_1) = \exp\left\{ - \frac{\delta}{2} \left[\alpha'\big(s(y_1,t_1)\big) - \alpha'\big(s(0,t_1)\big) \right] \right\}.
\label{eq:soln_A0}
\end{equation}
If $\alpha''(\mathfrak{t})$ is continuous, then, by the mean value theorem for integrals, we can write a more compact expression: $\mathcal{A}_0(y_1,t_1) = \exp\left\{-\tfrac{\delta}{2}\alpha''(\varsigma)y_1\right\}$, where $\varsigma \in \big( s(0,t_1), s(y_1,t_1) \big)$ is to be determined based on the functional form of $\alpha''$.

The ODE for the phase \eqref{eq:char_odes_psi0} gives $\overline \psi_0(s) = \overline\psi_0\big(s(0,t_1)\big) = 0$ $\forall t_1$ after applying the boundary condition \eqref{eq:ms_bc}. Substituting the positive solution from \eqref{eq:k_beta} for $k$, \eqref{eq:soln_A0} for $\mathcal{A}_0$ and $\psi_0=0$ into \eqref{eq:U0_solution} and then eliminating $s$ using \eqref{eq:s_t1_y1} completes the leading-order solution:
\begin{equation}
\hat{U}(\eta,\tau) \sim \exp\left\{ -\frac{\delta}{2}\Big[\alpha'\big( (1-\beta)\epsilon \tau - \delta\alpha(\epsilon \tau) + \epsilon \eta \big) - \alpha'\big( (1-\beta)\epsilon \tau - \delta\alpha(\epsilon \tau) \big) \Big] \right\}
\exp\left\{\ri\left[\tau - \frac{\eta}{1- \tilde\beta(\epsilon \tau)}\right]\right\}.
\label{eq:soln_ms_mv}
\end{equation}
Upon returning to the original (dimensional) variables in the stationary frame, we have
\begin{multline}
U(x,t) \sim \exp\Bigg\{ -\frac{\delta}{2}\Big[\alpha'\big( (1-\beta) \Omega t + x\Omega/c - \beta \Omega t - 2\delta \alpha(\Omega t) \big) - \alpha'\big( (1-\beta)\Omega t - \delta\alpha(\Omega t) \big) \Big] \Bigg\}\\
\times\exp\left\{\ri\frac{\omega}{1- \beta - \delta\alpha'(\Omega t)}\Big[t - x/c - \delta\alpha'(\Omega t)t + \delta\Omega^{-1} \alpha(\Omega t)\Big]\right\}.
\label{eq:soln_ms}
\end{multline}
This asymptotic result is valid for arbitrary $\delta$, provided that the speed of the emitter never exceeds the phase speed of waves in the medium [recall the discussion before \eqref{eq:var_param}], which corresponds to the requirement that the wave equation in the moving frame remains hyperbolic, namely $\delta < 1 - \beta$. Furthermore, we have to keep in mind that \eqref{eq:soln_ms} is valid only for $x < ct$ due to the finite speed of propagation of waves under \eqref{eq:wave_eq}. For $x > ct$, $U = 0$ since we stipulated homogeneous initial conditions.

Finally, we would like to make a brief comparison between \eqref{eq:soln_ms} and the corresponding solution  \cite[Eq.~(5.15)]{ow16} to the the IBVP for wave propagation in a homogeneous medium at rest, in which the acoustic source is modeled as a singular term on the right-hand side of \eqref{eq:wave_eq}. Although \eqref{eq:soln_ms} and \cite[Eq.~(5.15)]{ow16} look quite different, they have some conceptual similarities. The Doppler frequency including $\dot{X}_\mathrm{s}$ appears in the harmonic exponential term of both. Likewise, both \eqref{eq:soln_ms} and \cite[Eq.~(5.15)]{ow16} feature amplitude modulation but the functional form of the amplitude differs because \eqref{eq:soln_ms} is for a moving-boundary IBVP, while \cite[Eq.~(5.15)]{ow16} is for a moving source in a homogeneous medium at rest. Finally, while \cite[Eq.~(5.15)]{ow16} includes a summation over all solutions to the retarded time equation \cite[Eq.~(5.7)]{ow16}, \eqref{eq:soln_ms} does not because it is posed on a moving domain.


\section{Discussion}
\label{sec:discuss}

The prefactor of the bracketed expression inside the second exponential in the asymptotic solution \eqref{eq:soln_ms} gives the frequency of the observed wave at a distance $x$ from the emitter at time $t$. Thus, we learn that the emitted frequency $\omega$ is shifted to the observed frequency $\omega/[1-\beta - \delta \alpha'(\Omega t)]$, within $\mathcal{O}(\epsilon)$, due to the acceleration of the boundary. The shifted frequency has the same functional form as in the non-accelerating case \eqref{eq:doppler_shift}, except that $v/c \equiv \beta$ is replaced by $\dot X_\mathrm{s}(t)/c \equiv \beta + \delta\alpha'(\Omega t)$ [recall the second equation in \eqref{eq:wave_eq_bc}], consistent with \eqref{eq:doppler_shift_acc_heur}. Consequently, the wave experiences frequency \emph{modulation}\footnote{Frequency modulation in the Doppler spectrum of underwater acoustic waves due to source or receiver motion has been measured \cite{Keiff08}.} (rather than a simple shift) since $\alpha'$ depends upon $t$. In addition, within the same asymptotic order of approximation, there is an amplitude modulation of the waveform as embodied by first exponential in \eqref{eq:soln_ms}. Note, however, that all time dependences in \eqref{eq:soln_ms} are upon $\Omega t \equiv \epsilon \omega t$ [using \eqref{eq:var_param}], which is the slow time variable defined previously. Taking the limit $\delta\to0$ (no acceleration, $\dot X_\mathrm{s}(t)/c \to \beta = const.$), \eqref{eq:soln_ms_mv} becomes 
\begin{equation}
U(x,t) \sim \exp\left\{\ri\frac{\omega}{1- \beta}\Big[t - x/c\Big]\right\},
\end{equation}
which is a harmonic wave with constant amplitude and constant frequency, the latter given by the ordinary Doppler formula \eqref{eq:doppler_shift}, as required. 
Furthermore, note that for certain choices of $\alpha$, the first exponential in eqs.~\eqref{eq:soln_ms_mv} and \eqref{eq:soln_ms} could lead to the increase of the amplitude of the wave. 


To summarize: the effects on the observed waveform due to the acceleration of the boundary are (\textit{i}) amplitude modulation, as made explicit by the first exponential in eqs.~\eqref{eq:soln_ms_mv} and \eqref{eq:soln_ms}; (\textit{ii}) frequency modulation of the same functional form as for the case of an emitter with constant velocity; (\textit{iii}) a time-dependent phase shift, namely $\delta\alpha'(\Omega t)t - \delta\Omega^{-1}\alpha(\Omega t)$.

Next, based on \eqref{eq:soln_ms}, we can define two quantities that characterize the waveform: the frequency modulation (FM) factor
\begin{equation}
\mathsf{FM}(t) := \frac{1}{1- \beta - \delta \alpha'(\Omega t)}
\label{eq:FM_factor}
\end{equation}
and the amplitude modulation (AM) factor
\begin{equation}
\mathsf{AM}(x,t) := \exp\Bigg\{ -\frac{\delta}{2}\Big[\alpha'\big( (1-\beta)\Omega t + x\Omega/c - \beta \Omega t - 2\delta \alpha(\Omega t) \big) - \alpha'\big( (1-\beta)\Omega t - \delta\alpha(\Omega t) \big) \Big] \Bigg\}.
\end{equation}
Although our treatment applies to mechanical waves, we borrow this terminology from the radio wave transmission literature, wherein frequency modulation of the emitted wave at the source is used for ultra-high frequency (UHF) communications, while amplitude modulation of the emitted wave at the source is preferred for low frequency (LF) communications \cite{s97}. Additionally, as was the case with \eqref{eq:doppler_shift_acc_heur}, the FM factor \eqref{eq:FM_factor} is consistent with the FM factor found from the general three-dimensional solution for a moving source in a homogeneous stationary medium \cite[eq.~(5.23)]{ow16}.

$\mathsf{FM}$ and $\mathsf{AM}$ can be rewritten as functions of the three dimensionless parameters ($\beta$, $\delta$ and $\epsilon$), the dimensionless time $\tau=\omega t$ and the dimensionless distance $\omega x/c$:
\begin{subequations}
\begin{align}
\mathsf{FM}(t) &= \frac{1}{1- \beta - \delta \alpha'(\epsilon\omega t)}, \label{eq:FM_def}\\[2mm]
\mathsf{AM}(x,t) &= \exp\Bigg\{ -\frac{\delta}{2}\Big[\alpha'\big( \epsilon (1-\beta)\omega t + \epsilon \omega x/c - \epsilon \beta \omega t - 2\delta \alpha(\epsilon\omega t) \big)  \label{eq:AM_def} - \alpha'\big( (1-\beta)\epsilon\omega t - \delta\alpha(\epsilon\omega t) \big) \Big] \Bigg\}.
\end{align}\label{eq:AM_FM_def}\end{subequations}
Note that while the frequency modulation is only a function of time, the amplitude modulation is both a function of time and space. Consequently, the wave amplitude measured by an observer depends on the observer's instantaneous distance from the emitter.

Next, we consider two illustrative examples.

\subsection{A decelerating boundary}
\label{sec:decel}

First, consider the case of a continuously decelerating boundary 
We take
\begin{equation}
\alpha(\mathfrak{t}) = 1 - \re^{-\mathfrak{t}}
\label{eq:at_decel}
\end{equation}
so that $\alpha'(\mathfrak{t}) = \re^{-\mathfrak{t}}$, $\max_{\mathfrak{t}\ge0} |\alpha'(\mathfrak{t})|=1$ and $\alpha(0) = 0$. Then, evaluating \eqref{eq:AM_FM_def} using \eqref{eq:at_decel} yields
\begin{subequations}\begin{align}
\mathsf{FM}(t) &= \frac{1}{1- \beta - \delta \re^{-\epsilon\omega t}},\\[2mm]
\mathsf{AM}(x,t) &= \exp\left\{\frac{\delta}{2} \re^{\delta - 2 \delta \re^{-\epsilon\omega t}} \left[\re^{\delta\re^{-\epsilon\omega t} - \epsilon (1 - \beta)\omega t} - \re^{\delta - \epsilon\omega t - \epsilon\omega x/c + 2\epsilon\beta\omega t}\right]\right\}.\label{eq:AM_decel}
\end{align}\label{eq:AM_FM_decel}\end{subequations}
The $\mathsf{FM}$ factor is not singular thanks to our earlier restriction of $\delta < 1-\beta$.

We can compute the following asymptotic limits of \eqref{eq:AM_FM_decel}:
\begin{subequations}\begin{align}
\mathsf{FM} &\to \begin{cases} \displaystyle\frac{1}{1 - \beta - \delta}, &t \to 0;\\[4mm]
\displaystyle\frac{1}{1 - \beta}, &t \to \infty;\end{cases}\\[2mm]
\mathsf{AM} &\to \begin{cases} \displaystyle \exp\left\{\frac{\delta}{2}(1-\re^{-\epsilon\omega x/c})\right\}, &t \to 0;\\[4mm]
1, &t \to \infty.\end{cases}
\end{align}\label{eq:AM_FM_limits_decel}\end{subequations}
As $t\to\infty$, we observe that $\mathsf{FM}$ and $\mathsf{AM}$ reduce to the classical Doppler  relations. This is expected as an exponentially-decaying acceleration quickly becomes negligible.

\begin{figure}
	\centering
	\includegraphics[width=0.75\textwidth]{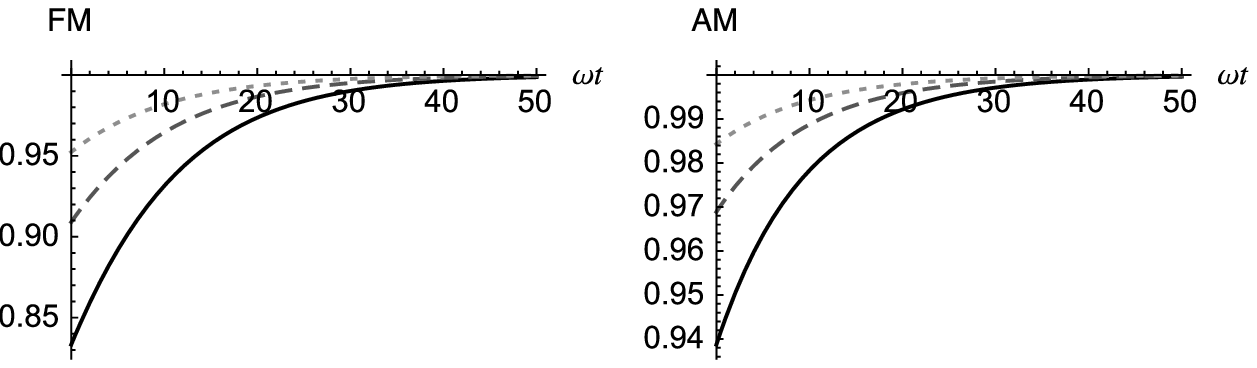}
	\caption{Frequency (left) and amplitude (right, at $\omega x/c = 10$) modulation factors for $\epsilon = 0.1$ and $\beta=0$; $\delta=-0.2$ (solid), $\delta=-0.1$ (dashed) and $\delta=-0.05$ (dotted).} 
	\label{fig:decel}
\end{figure}

\begin{figure}
	\centering
	\includegraphics[width=0.75\textwidth]{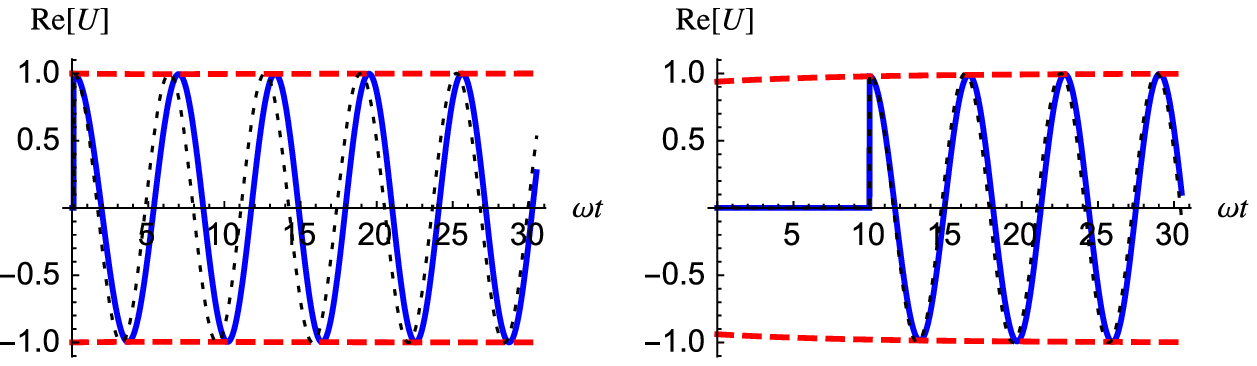}
	\caption{Waveform $\real[\mathsf{AM}(x,t)\re^{\ri\mathsf{FM}(t)\omega(t-x/c)}]$ observed (solid curves) by two receivers in the fixed frame, as a function of dimensionless time; $\epsilon=0.1$, $\beta = 0$, $\delta = -0.2$; $\omega x/c = 0.1$ (left) and $\omega x/c = 10$ (right). Dashed curves represent the envelope $\mathsf{AM}$ from \eqref{eq:AM_decel}, while the dotted curves represent the harmonic wave $\real[\re^{\ri \omega (t - x/c)}]$, as it would propagate away from a stationary source.}
	\label{fig:ReU_decel}
\end{figure}

Figure~\ref{fig:decel} shows plots of $\mathsf{FM}$ and $\mathsf{AM}$ for some representative values of the dimensionless parameters. Meanwhile, fig.~\ref{fig:ReU_decel} shows the waveform $\real[\mathsf{AM}(x,t)\re^{\ri\mathsf{FM}(t)\omega(t-x/c)}]$ at fixed $\omega x/c$. To simplify the discussion, the time-dependent phase shift has been neglected and, without loss of generality, $\beta = 0$ is used in these figures. The plots show that the effects due to acceleration disappear as $t\to\infty$ (\textit{i.e.}, as the boundary velocity becomes uniform or, in this case of $\beta=0$, the boundary becomes stationary) as evidenced by the increasing overlap between the solid curves, $\real[\mathsf{AM}(x,t)\re^{\ri\mathsf{FM}(t)\omega(t-x/c)}]$, and dotted curves, $\real[\re^{\ri\omega(t-x/c)}]$, in fig.~\ref{fig:ReU_decel} for large $\omega t$. This observation is supported by the long-time asymptotics given by \eqref{eq:AM_FM_limits_decel}. As can be seen in fig.~\ref{fig:ReU_decel}, the frequency shift due to the acceleration of the boundary is more pronounced for $\omega x/c \ll 1$ (close to the source), while the amplitude shift is more pronounced for $\omega x/c \gg 1$ (farther downstream).

\begin{figure}
	\centering
	\includegraphics[width=0.75\textwidth]{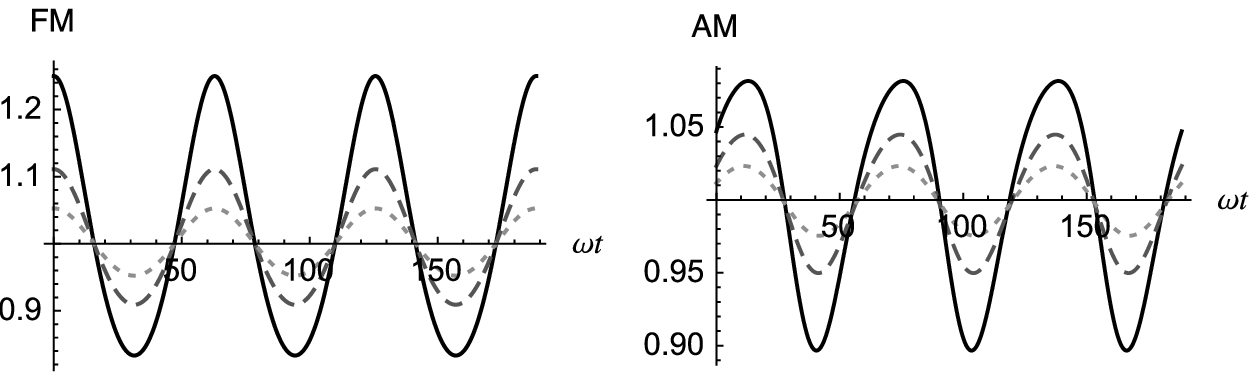}
	\caption{Frequency (left) and amplitude (right, at $\omega x/c = 10$) modulation factors for $\epsilon = 0.1$ and $\beta=0$; $\delta=0.2$ (solid), $\delta=0.1$ (dashed) and $\delta=0.05$ (dotted).}
	\label{fig:osc}
\end{figure}

\begin{figure}
	\centering
	\includegraphics[width=0.75\textwidth]{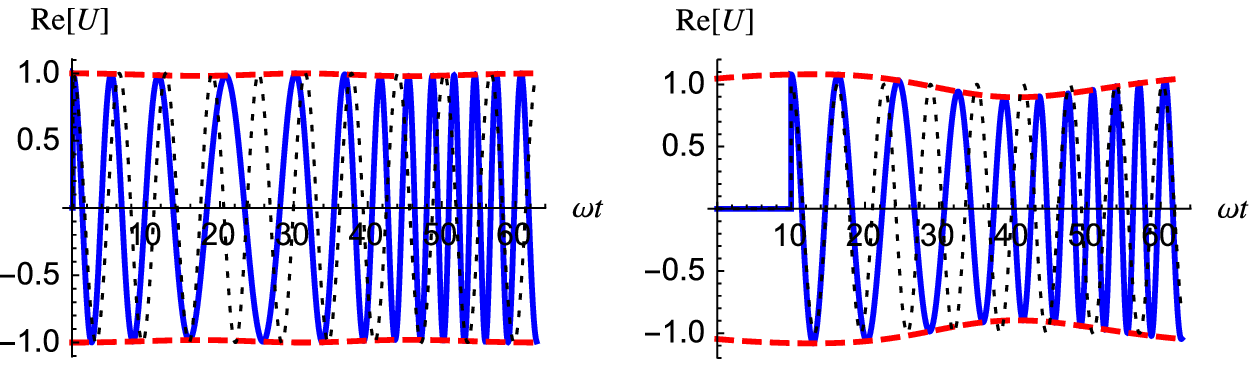}
	\caption{Waveform $\real[\mathsf{AM}(x,t)\re^{\ri\mathsf{FM}(t)\omega(t-x/c)}]$ observed (solid curves) by two receivers in the fixed frame, as a function of dimensionless time; $\epsilon=0.1$, $\beta = 0$, $\delta = 0.2$; $\omega x/c = 0.1$ (left) and $\omega x/c = 10$ (right). Dashed curves represent the envelope $\mathsf{AM}$ from \eqref{eq:AM_osc}, while the dotted curves represent the harmonic wave $\real[\re^{\ri \omega (t - x/c)}]$, as it would propagate away from a stationary source.}
	\label{fig:ReU_osc}
\end{figure}

\subsection{Periodic oscillations of the boundary}
\label{sec:osc}

Second, consider an oscillating boundary motion.
In this case, we take
\begin{equation}
\alpha(\mathfrak{t}) = \sin(\mathfrak{t})
\label{eq:at_osc}
\end{equation}
so that $\alpha'(\mathfrak{t}) = \cos(\mathfrak{t})$, $\max_{\mathfrak{t}\ge0}|\alpha'(\mathfrak{t})| = 1$ and $\alpha(0)=0$. Evaluating \eqref{eq:AM_FM_def} using \eqref{eq:at_osc} yields
\begin{subequations}\begin{align}
\mathsf{FM}(t) &= \frac{1}{1-\beta-\delta \cos(\epsilon\omega t)},\label{eq:FM_osc}\\[2mm]
\mathsf{AM}(x,t) &= \exp\Big\{-\frac{\delta}{2} \Big[\cos\big(\epsilon(1-\beta) \omega t + \epsilon \omega x/c - \epsilon \beta \omega t - 2 \delta \sin(\epsilon \omega t)\big)\label{eq:AM_osc}
- \cos\big(\epsilon(1 - \beta) \omega t - \delta \sin(\epsilon\omega t)\big)\Big]\Big\}.
\end{align}\end{subequations}
Once again, $\delta < 1 - \beta$ guarantees that the latter expressions are free of singularities. The expression \eqref{eq:FM_osc} for $\mathsf{FM}$ agrees with the heuristic expression in \cite[eq.~(6)]{a07}. These time-periodic frequency and amplitude modulations can also be loosely interpreted as combined vibrato and tremolo produced by, \textit{e.g.,} a Leslie loudspeaker \cite[pp.~521--523]{CG87}.

The standard Doppler effect, as quantified by $\beta$, is to change the mean of the $\mathsf{FM}$ and $\mathsf{AM}$ factors, so we take $\beta=0$ for the remainder of this subsection without loss of generality. This assumption corresponds to an oscillating boundary with zero net displacement. As in \S\ref{sec:decel}, figures~\ref{fig:osc} and \ref{fig:ReU_osc} illustrate, respectively, the frequency and amplitude modulation factors, \textit{i.e.}, $\mathsf{FM}$ and $\mathsf{AM}$, and the waveform $\real[\mathsf{AM}(x,t)\re^{\ri\mathsf{FM}(t)\omega(t-x/c)}]$ at fixed $\omega x/c$, neglecting the time dependent phase shift and with $\beta = 0$. Once again, fig.~\ref{fig:ReU_osc} highlights the fact that the amplitude shift due to the acceleration of the boundary is more pronounced for $\omega x/c \gg 1$ (far downstream from the source). In this example, however, the frequency modulation remains a pernicious feature for all $t$ because $\mathsf{FM}(t)\nrightarrow 1/(1-\beta) = 1$ (for $\beta =0$) as $t\to\infty$.


\section{Conclusion}
\label{sec:conclusion}

In this work, we investigated the effect of the acceleration of a moving boundary of a 1D domain on the emitted mechanical plane waves from a co-located source in a homogeneous resting medium, which represents a new variant of George Carrier's ``spaghetti problem'' \cite{c49}. The hyperbolic partial differential  equation describing the wave motion was posed as a boundary-value problem with a moving ``inlet'' boundary condition at the emitter's position (\textit{i.e.}, at the moving boundary) and a radiation condition at infinity. This problem was transformed to a dispersive hyperbolic PDE with \emph{variable} coefficients in the moving frame.

The small parameter $\epsilon := \Omega/\omega$, which represents the ratio of the characteristic time scale of the emitted wave to the characteristic time scale on which the mechanical oscillations of the boundary take place, was used in a multiple-scales asymptotic expansion of the solution of the boundary-value problem in the moving frame. Specifically, we derived the solution given in \eqref{eq:soln_ms}, which consists of an envelope propagating over the Doppler shifted carrier wave. This solution is a {formal} new result for moving IBVPs. 
Futhermore, since the  governing PDE in the moving frame \eqref{eq:wave_eq_mf} has variable coefficients, it also has, in particular, a variable phase speed; thus, the results in \S\ref{sec:ms} generalize the problem and solution from the appendix of \cite{RabbittHolmes}. 

It is interesting to note that, in the case of a decelerating boundary motion considered in \S\ref{sec:decel}, the amplitude modulation is pure attenuation. A similar a effect is experienced by a harmonic wave traveling through a porous medium under the nonlinear theory of acoustics \cite{j09}. In the latter case the source is stationary, but the medium through which the wave travels induces signal loss, and the signal itself experiences nonlinear effects due to compressibility. Of course, in the present work, the attenuation is purely due to the motion of the boundary (it arises in the absence of dissipative or nonlinear effects, although it would also be of interest to consider the case of complex dispersive moving media \cite{gvn09}).

More generally, our discussion of the Doppler effect in the context of moving IBVPs has important implications for a number of other wave phenomena. Inhomogeneities of the carrier medium can cause diffraction of acoustic and/or elastic waves, which can be interrogated using techniques similar to the present mathematical framework \cite{kk59}. 
In plasmas, the variable refractive index of the medium can also cause frequency modulations that are not explainable by the classical Doppler effect \cite{dss00}.

We note that, in our illustrated examples in \S\ref{sec:discuss}, we neglected the effect of the time-dependent phase shift caused by the acceleration of the boundary. It is conceivable that this feature could explain the experimental observations \cite{bf01,a07} of generation of higher harmonics in the spectrum of the received signal. Future work also includes studying the evolution of a Gaussian wave packet under \eqref{eq:wave_eq_mf} in the spirit of \cite{c08}. 

Finally, we emphasize that our work does \emph{not} treat the relativistic case, \textit{i.e.}, the field of moving electromagnetic charges, which has been examined in detail in textbooks \cite[Chap.~8]{LL80}.

\acks
I.C.C.\ was partially supported by the LANL/LDRD Program through a Feynman Distinguished Fellowship. I.C.C.\ thanks an anonymous reviewer for incisive comments, which significantly improved the presentation, and for providing key references.

\bibliographystyle{wileyj}
\bibliography{accel_doppler}


\section*{Appendix. Formulation using characteristic coordinates}
\setcounter{equation}{0}
\renewcommand{\theequation}{{A.\arabic{equation}}}

Consider again the governing PDE \eqref{eq:wave_eq_nd}:
\begin{equation}
\frac{\partial^2 \hat{U}}{\partial \tau^2} - 2\tilde\beta(\tau_1)\frac{\partial^2 \hat{U}}{\partial\eta\partial\tau} 
 - \big[1 - \tilde\beta^2(\tau_1) \big] \frac{\partial^2 \hat{U}}{\partial \eta^2} 
 - \epsilon\delta\alpha''(\tau_1)\frac{\partial \hat{U}}{\partial\eta} = 0,\qquad
 0<\eta<\infty,\quad 0<\tau<\infty.
\label{eq:wave_eq_nd_2}
\end{equation}
In general, we can impose ``initial'' conditions of the form
\begin{equation}
\hat{U}(0,\tau) = \mathring{U}(\tau), \qquad \DD{\hat{U}}{\eta}(0,\tau) = 0,
\label{eq:ic}
\end{equation}
where $\mathring{U}(\cdot)$ is some given excitation, and $\eta$ plays the role of the ``time-like" variable.


Without the simplifying assumption of a harmonic excitation, implementing the multiple-scale expansion proceeds by introducing characteristic coordinates for the slow scales only (see, \textit{e.g.}, \cite[\S3.9]{h95} and also \cite{hv96} for illuminating applications of this technique):
\begin{equation}
\theta_1 = \eta + (\tilde \beta+1) \tau, \qquad \theta_2 = \eta + (\tilde \beta -1) \tau, \qquad \tau_1 = \epsilon \tau, \qquad \eta_1 = \epsilon \eta.
\end{equation}
As before, replacing derivates using the chain rule and letting $\hat{U}(\eta,\tau) = \mathcal{U}(\theta_1,\theta_2,\eta_1,\tau_1)$ transforms \eqref{eq:wave_eq_nd_2} is into
\begin{equation}
-4 \DDM{\mathcal{U}}{\theta_1}{\theta_2} = -2\epsilon \DD{}{\tau_1} \left( \DD{\mathcal{U}}{\theta_1} - \DD{\mathcal{U}}{\theta_2} \right)
+ 2\epsilon \DD{}{\eta_1} \left[(1 + \tilde\beta)\DD{\mathcal{U}}{\theta_1} + (1 - \tilde\beta)\DD{\mathcal{U}}{\theta_2} \right]
+ \epsilon\delta\alpha''\left(\DD{\mathcal{U}}{\theta_1} + \DD{\mathcal{U}}{\theta_2}\right) + \mathcal{O}(\epsilon^2).
 \label{eq:characteristics_eps}
\end{equation}
The boundary condition is specified at $\eta=0$, which corresponds to $\theta_1 = (\tilde\beta+1)\tau$, $\theta_2 = (\tilde\beta-1)\tau$ and $\eta_1=0$. Therefore, \eqref{eq:ic} becomes
\begin{equation}
\mathcal{U}(\theta_1,\theta_2,0,\tau_1) = \mathring{U}\left(\frac{\theta_1}{\tilde\beta+1}\right) = \mathring{U}\left(\frac{\theta_2}{\tilde\beta -1}\right),\qquad
\left(\DD{}{\theta_1} + \DD{}{\theta_2}+ \epsilon\DD{}{\eta_1} + \cdots\right)\mathcal{U}(\theta_1,\theta_2,0,\tau_1) = 0.
\label{eq:ic_ms}
\end{equation}

We proceed by a regular expansion of the independent variable
\begin{equation}
\mathcal{U}(\theta_1,\theta_2,\eta_1,\tau_1) = \mathcal{U}_0(\theta_1,\theta_2,\eta_1,\tau_1) + \epsilon\, \mathcal{U}_1(\theta_1,\theta_2, \eta_1,\tau_1) + \cdots.
\end{equation}
At the leading order, we obtain the following ``initial-value'' problem:
\begin{subequations}\begin{align}
&\DDM{\mathcal{U}_0}{\theta_1}{\theta_2} = 0, \label{eq:th1_th2_ivp_0_pde}\\
&\mathcal{U}_0 (\theta_1,\theta_2,0,\tau_1) = \mathring{U}\left(\frac{\theta_1}{\tilde\beta+1}\right), \qquad \left(\DD{}{\theta_1} + \DD{}{\theta_2}\right)\mathcal{U}_0 (\theta_1,\theta_2,0,\tau_1) =0. \label{eq:th1_th2_ivp_0_ic}
\end{align}\label{eq:th1_th2_ivp_0}\end{subequations} 
The general d'Alembert-type solution to \eqref{eq:th1_th2_ivp_0} is
\begin{equation}
\mathcal{U}_0 = \mathcal{G}_{0,1}(\theta_1,\eta_1,\tau_1) + \mathcal{G}_{0,2}(\theta_2,\eta_1,\tau_1),
\label{eq:leading_order_2}
\end{equation}
where $\mathcal{G}_{0,1}$ and $\mathcal{G}_{0,2}$ are to be determined.

The $\mathcal{O}(\epsilon)$ problem is
\begin{equation}
-4 \DDM{\mathcal{U}_1}{\theta_1}{\theta_2} = -2 \DD{}{\tau_1} \left( \DD{\mathcal{U}_0}{\theta_1} - \DD{\mathcal{U}_0}{\theta_2} \right)
+ 2 \DD{}{\eta_1} \left[(1 + \tilde\beta)\DD{\mathcal{U}_0}{\theta_1} + (1 - \tilde\beta)\DD{\mathcal{U}_0}{\theta_2} \right]
+ \delta\alpha''\left(\DD{\mathcal{U}_0}{\theta_1} + \DD{\mathcal{U}_0}{\theta_2}\right).
\label{eq:app_oeps_pb}
\end{equation}
Substituting the $\mathcal{O}(1)$ solution from \eqref{eq:leading_order_2} into \eqref{eq:app_oeps_pb} and solving, we find that
\begin{multline}
\mathcal{U}_1 = \mathcal{G}_{1,1}(\theta_1,\eta_1,\tau_1) + \mathcal{G}_{1,2}(\theta_2,\eta_1,\tau_1)\\
- \frac{1}{4}\theta_2\left(-2 \DD{}{\tau_1} + 2 (1 + \tilde\beta)\DD{}{\eta_1} + \delta\alpha'' \right) \mathcal{G}_{0,1}
- \frac{1}{4}\theta_1\left(2 \DD{}{\tau_1} + 2 (1 - \tilde\beta)\DD{}{\eta_1} + \delta\alpha'' \right) \mathcal{G}_{0,2},
\end{multline}
where $\mathcal{G}_{1,\{1,2\}}(\theta_{1,2},\eta_1,\tau_1)$ satisfy the homogeneous PDE for $\mathcal{U}_1$.
To suppress secular terms, we must require that
\begin{subequations}\begin{align}
-\DD{\mathcal{G}_{0,1}}{\eta_1} + \frac{1}{1 + \tilde\beta(\tau_1)} \DD{\mathcal{G}_{0,1}}{\tau_1} &= - \frac{\delta\alpha''(\tau_1)}{2[1 + \tilde\beta(\tau_1)]} \mathcal{G}_{0,1},\\
\DD{\mathcal{G}_{0,2}}{\eta_1} + \frac{1}{1 - \tilde\beta(\tau_1)}\DD{\mathcal{G}_{0,2}}{\tau_1} &= - \frac{\delta\alpha''(\tau_1)}{2[1 - \tilde\beta(\tau_1)]} \mathcal{G}_{0,2}.
\end{align}\end{subequations}
This is a pair of uncoupled scalar hyperbolic PDEs. We proceed by the method of characteristics:
\begin{equation}
\frac{\rd \eta_1}{\rd s_{1,2}} = \mp 1,\qquad \frac{\rd \tau_1}{\rd s_{1,2}} = \frac{1}{1 \pm \tilde\beta(\tau_1)},\qquad \frac{\rd \mathcal{G}_{0,\{1,2\}}}{\rd s_{1,2}} =  - \frac{\delta\alpha''(\tau_1)}{2[1 \pm \tilde\beta(\tau_1)]} \mathcal{G}_{0,\{1,2\}}.
\label{eq:characs_tau1_eta1}
\end{equation}
Recalling that $\tilde\beta(t_1) \equiv \beta + \delta \alpha'(t_1)$, the first two equations in \eqref{eq:characs_tau1_eta1} give
\begin{equation}
s_{1,2}(\eta_1,\tau_1) = \int_0^{\tau_1} 1\pm \tilde\beta(\mathfrak{t}) \,\rd\mathfrak{t} \mp\eta_1 = (1\pm\beta)\tau_1 \pm \delta\alpha(\tau_1) \mp \eta_1.
\label{eq:s_eta1_tau1}
\end{equation}
Integrating the third ODE in \eqref{eq:characs_tau1_eta1}, we obtain
\begin{equation}
\int_{s_{1,2}(0,\tau_1)}^{s_{1,2}(\eta_1,\tau_1)} \frac{\rd \mathcal{G}_{0,\{1,2\}}}{\mathcal{G}_{0,\{1,2\}}} = - \int_{s_{1,2}(0,\tau_1)}^{s_{1,2}(\eta_1,\tau_1)}  \frac{\delta\alpha''(\tau_1(\varsigma))}{2[1 \pm \tilde\beta(\tau_1(\varsigma))]} \,\rd \varsigma = - \frac{\delta}{2} \int_{s_{1,2}(0,\tau_1)}^{s_{1,2}(\eta_1,\tau_1)} \alpha''(\tau_1) \,\rd \tau_1,
\label{eq:app_charac_int}
\end{equation}
where we used the second equation in \eqref{eq:characs_tau1_eta1} to change variables in the integral and the lower limit was chosen because the initial condition \eqref{eq:ic_ms} is specified at $\eta_1=0$.
Finally, performing the integration in \eqref{eq:app_charac_int} and eliminating $s_{1,2}$ using \eqref{eq:s_eta1_tau1}, we arrive at
\begin{equation}
\mathcal{G}_{0,\{1,2\}}(\theta_{1,2},\eta_1,\tau_1) = \mathring{\mathcal{G}}_{0,\{1,2\}}\big(\theta_{1,2},\tau_1\big)
\exp\Bigg\{ -\frac{\delta}{2}\big[\alpha'\big((1\pm\beta)\tau_1 \pm \delta\alpha(\tau_1) \mp \eta_1\big) - \alpha'\big((1\pm\beta)\tau_1 \pm \delta\alpha(\tau_1)\big)\big] \Bigg\},
\end{equation}
where $\mathring{\mathcal{G}}_{0,\{1,2\}}\big(\theta_{1,2},\tau_1\big)$ are uniquely determined by the ``initial'' condition at $\eta_1=0$, namely \eqref{eq:th1_th2_ivp_0_ic}.
The amplitude modulation is identical to that in \eqref{eq:soln_ms_mv}, while the frequency modulation (and classical Doppler shift) are already ``built into'' the characteristic variables $\theta_{1,2}$.

%

\end{document}